\documentclass[letterpaper]{amsart} 

\usepackage{amsmath}
\usepackage[all]{xy}

\usepackage{psfrag}
\usepackage{epsfig}
 
\numberwithin{equation}{section}

\newtheorem{theorem}{Theorem}[section]

\newtheorem{corollary}[theorem]{Corollary}

\theoremstyle{definition}

\def\cc{\mathbf{c}}

\def\TT{\mathbb{T}}

\def\ZZ{\mathbb{Z}}

\newcommand{\overunder}[2]{
\!\begin{array}{c}
\scriptstyle{#1}\\[-.1in]
-\!\!\!-\!\!\!-\\[-.1in]
\scriptstyle{#2}
\end{array}
\!
}
\begin{document}

\title{Cluster algebras and symmetrizable matrices}

\author{Ahmet I. Seven}

\address{Middle East Technical University, Department of Mathematics, 06800, Ankara, Turkey}
\email{aseven@metu.edu.tr}

\thanks{The author's research was supported in part by the Turkish Research Council (TUBITAK) }


\date{February 14, 2018}



\begin{abstract}

In the structure theory of cluster algebras, principle coefficients are parametrized by a family of integer vectors, called {$\cc$-vectors}. Each $\cc$-vector with respect to an acyclic initial seed is a real root of the corresponding root system and the $\cc$-vectors associated with any seed defines a symmetrizable quasi-cartan companion for the corresponding exchange matrix. We establish basic combinatorial properties of these companions. 
In particular, we show that $\cc$-vectors define an admissible cut of edges in the associated diagrams.


\end{abstract}


\maketitle

\section{Introduction}
\label{sec:intro}


In the structure theory of cluster algebras, principle coefficients are parametrized by a family of integer vectors, called {$\cc$-vectors}. Each $\cc$-vector with respect to an acyclic initial seed is a real root of the corresponding root system; furthermore, the $\cc$-vectors associated with any seed defines a symmetrizable quasi-cartan companion for the corresponding exchange matrix \cite[Corollary~3.29]{RS}. In this paper, we study basic combinatorial properties of these companions. In particular, we show that $\cc$-vectors define an admissible cut of edges in the associated diagrams.


To state our results, we need some terminology. Let us recall that an $n\times n$ integer matrix $B$ is skew-symmetrizable if there is a diagonal matrix $D$ with positive diagonal entries such that $DB$ is skew-symmetric. 
We denote by $\TT_n$ an \emph{$n$-regular tree} whose edges are labeled by the numbers $1, \dots, n$ such that the $n$ edges incident to each vertex have different labels. The notation $t \overunder{k}{} t'$ indicates that vertices $t,t'\in\TT_n$ are connected by an edge labeled by~$k$. We fix a vertex $t_0$ in $\TT_n$ and assign the pair $(\cc_0,B_0)$, where $\cc_0$ is the tuple of standard basis and $ B_0 $ is a skew-symmetrizable matrix. Then,  to every vertex $t \in \TT_n$ we assign a pair,  called  a $Y$-seed, $(\cc_t, B_t)$, where $\cc_t=(\cc_1,...,\cc_n)$ with each $\cc_i=\cc_{i;t}=(c_1,...,c_n) \in \ZZ^n$ being non-zero and having either all entries nonnegative or all entries nonpositive; we write $sgn(\cc_i)=+1$ or $sgn(\cc_i)=-1$ respectively and call it a \emph{$\cc$-vector}. Furthermore, for any edge  $t \overunder{k}{} t'$,  the $Y$-seed $(\cc', B')=(\cc_{t'}, B_{t'})$ is obtained from 
$(\cc, B)=(\cc_{t}, B_{t})$ by the \emph{$Y$-seed mutation} $\mu_k$ defined as follows, where we denote $[b]_+ = \max(b,0)$:
\begin{itemize}
\item
The entries of the matrix $B'=(B'_{ij})$ are given by
\begin{equation}
\label{eq:matrix-mutation}
B'_{ij} =
\begin{cases}
-B_{ij} & \text{if $i=k$ or $j=k$;} \\[.05in]
B_{ij} + [B_{ik}]_+ [B_{kj}]_+ - [-B_{ik}]_+ [-B_{kj}]_+
 & \text{otherwise.}
\end{cases}
\end{equation}
\item
The tuple $\cc'=(\cc_1',\dots,\cc_n')$ is given by
\begin{equation}
\label{eq:y-mutation}
\cc'_i =
\begin{cases}
-\cc_{i} & \text{if $i = k$};\\[.05in]
\cc_i+[sgn(\cc_k)B_{k,i}]_+\cc_k
 & \text{if $i \neq k$}.
\end{cases}
\end{equation}
\end{itemize}
By \cite[Corollary 5.5]{GHKK}, each $\cc'_i=(c'_1,...,c'_n) $ also has either all entries nonnegative or all entries nonpositive. The matrix $B'$ is skew-symmetrizable with the same choice of $D$; we write $B' = \mu_k(B)$ and call the transformation $B \mapsto B'$ the \emph{matrix mutation}. 
For the $Y$-seeds, we denote $\mu_k(\cc, B)=(\cc', B')$; we call $(\cc_0,B_0)$ the initial $ Y $-seed. It is well known that mutation is an involutive operation.

Let us also recall that
the \emph{diagram} of a skew-symmetrizable $n\times n$ matrix ${B}$ is the directed graph $\Gamma ({B})$ whose vertices are the indices $1,2,...,n$ such that there is a directed edge from $i$ to $j$ if and only if ${B}_{j,i} > 0$, and this edge is assigned the weight $|B_{ij}B_{ji}|\,$. The diagram $\Gamma(B)$ is called acyclic if it has no oriented cycles. Then there is a corresponding generalized Cartan matrix $A$ such that $A_{i,i}=2$ and $A_{i,j}=-|B_{i,j}|$ for $i\ne j$. There is also the associated root system in the root lattice spanned by the simple roots $\alpha_i$ \cite{K}. For each simple root $\alpha_i$, the corresponding reflection $s_{\alpha_i}=s_i$ is the linear isomorphism defined on the basis of simple roots as $s_i(\alpha_j)=\alpha_j-A_{i,j}\alpha_i$. Then the real roots are defined as the vectors obtained from the simple roots by a sequence of reflections. It is well known that the coordinates of a real root with respect to the basis of simple roots are either all nonnegative or all nonpositive, see \cite{K} for details.

On the other hand, an $n\times n$  matrix $A$ is called symmetrizable if there exists a symmetrizing diagonal matrix $D$ with positive diagonal entries such that $DA$ is symmetric. 
A \emph{quasi-Cartan companion} (or "companion" for short) of a skew-symmetrizable matrix $B$ is a symmetrizable matrix $A$ such that $A_{i,i}=2$, $|A_{i,j}|= |B_{i,j}|$ for all $i \ne j$. 

A fundamental relation between $Y$-seeds and symmetrizable matrices has been given in  \cite[Corollary~3.29]{RS} as follows: 

\begin{theorem}\label{th:DWZ2-RS} \cite[Corollary~3.29]{RS}  Suppose that the initial seed $(\cc_0,B_0)$ is acyclic. Then, for any $Y$-seed $ (\cc_t, B_t)$, $t \in \TT_n$, each $\cc$-vector $\cc_{i}=\cc_{i;t}$ is the coordinate vector of a real root with respect to the basis of simple roots in the corresponding root system. Furthermore, $A=A_t=(\langle\cc_j,{\cc}^{\vee}_i\rangle)$,  the matrix of the pairings between the roots and the coroots, is a quasi-Cartan companion of the skew-symmetrizable matrix $ B=B_t $. 

(The matrices $ A_t $ are symmetrizable with the same choice of a symmetrizing matrix $D$, which is also skew-symmetrizing for all $ B_t $.)
\end{theorem}

An important combinatorial property related to quasi-Cartan companions is admissibility \cite{S3, S6}, which is a generalization of the notion of a generalized Cartan matrix. More precisely, a quasi-Cartan companion $A$ of a skew-symmetrizable matrix $B$ \emph{admissible} if, for any oriented (resp. non-oriented) cycle $Z$ in $\Gamma(B)$, there is exactly an odd (resp. even) number of edges $\{i,j\}$ such that $A_{i,j}>0$.  If $\Gamma(B)$ is acyclic, then the associated  generalized Cartan matrix is admissible. Our first result generalizes this property by showing that the quasi-Cartan companions defined by $\cc$-vectors are also admissible:


\begin{theorem}\label{th:companion2}
In the set-up of Theorem \ref{th:DWZ2-RS}, the quasi-Cartan companion $A$  has the following properties:

\begin{itemize}
\item
Every directed path of the diagram $\Gamma(B)$ has at most one edge $\{i,j\}$ such that $A_{i,j}>0$.
\item
Every oriented cycle of the diagram $\Gamma(B)$ has exactly one edge $\{i,j\}$ such that $A_{i,j}>0$.

\item
Every non-oriented cycle of the diagram $\Gamma(B)$ has an even number of edges $\{i,j\}$ such that $A_{i,j}>0$.

\end{itemize}
In particular, the quasi-Cartan companion $A$ is admissible. Furthermore, any admissible quasi-Cartan companion of $B$ can be obtained from $A$ by a sequence of simultaneous sign changes in rows and columns. 

\end{theorem}
\noindent
The special case of this theorem when $ B $ is skew-symmetric was obtained in \cite[Theorem 1.4]{S6} by the author.  Let us also recall from \cite{S6} that a set $C$ of edges in $\Gamma(B)$ is called an ''admissible cut'' if every oriented cycle contains exactly one edge that belongs to $C$ and every non-oriented cycle contains exactly an even number of edges in $C$. Thus, in the setup of the theorem, the $\cc$-vectors define an admissible cut of edges: the set of edges $\{i,j\}$ in $\Gamma(B)$ such that $A_{i,j}>0$ is an admissible cut. For skew-symmetric matrices, this notion has been applied to the representation theory of algebras in \cite{HI,BRS}. 

Our next result is the following explicit description of the quasi-Cartan companions defined by the $\cc$-vectors:





\begin{theorem}\label{th:companion0}
In the set-up of Theorem \ref{th:DWZ2-RS}, the quasi-Cartan companion $A$  has the following properties:
\begin{itemize}
\item
If $sgn(B_{j,i})=sgn(\cc_j)$, then $A_{j,i}=-sgn(\cc_j)B_{j,i}=-|B_{j,i}| $.
\item
If $sgn(B_{j,i})=-sgn(\cc_j)$, then $A_{j,i}=sgn(\cc_i)B_{j,i}=-sgn(\cc_i)sgn(\cc_j)|B_{j,i}|$.
\end{itemize}

\smallskip
In particular; if $sgn(\cc_j)=-sgn(\cc_i)$, then $B_{j,i}=sgn(\cc_i)A_{j,i}$.







\end{theorem}

\noindent
Let us note that the special case of this theorem when $ B $ is skew-symmetric was obtained in \cite[Theorem 1.3]{S6} by the author. We will prove this more general theorem 
using \cite[Corollary~3.29]{RS}, which has been given above as Theorem \ref{th:DWZ2-RS}. (Note that the statement \cite[Corollary~3.29]{RS} was not present in the earlier versions of \cite{RS}).

\begin{corollary}\label{cor:companion}
In the setup of Theorem \ref{th:companion0}, suppose that $t \overunder{k}{} t'$ in $\TT_n$ 
Then, for $\mu_k(\cc, B)=(\cc', B')$, we have the following: 
if $\cc'_i\ne \cc_i$, then $\cc'_i=s_{\cc_k}(\cc_i)$, where $s_{\cc_k}$ is the reflection with respect to the real root $\cc_k$ and $\ZZ^n$ is identified with the root lattice.



\end{corollary}

\begin{figure}[ht]

\setlength{\unitlength}{2.6pt}

\begin{center}

\begin{picture}(60,30)(-30,0)

\put(-32,15){\makebox(0,0){$2$}}
\put(-23,20){\makebox(0,0){$2$}}
\put(-13,25){\makebox(0,0){$2$}}

\thicklines

\put(-30,0){\circle*{2.0}}
\put(-15,15){\circle*{2.0}}
\put(-30,30){\circle*{2.0}}
\put(0,15){\circle*{2.0}}




\put(-30,0){\vector(0,1){30}}
\put(-30,30){\vector(2,-1){30}}
\put(-30,30){\vector(1,-1){15}}
\put(-30,0){\vector(2,1){30}}
\put(-15,15){\vector(-1,-1){15}}
\put(-15,15){\vector(1,0){15}}

\end{picture}

\end{center}

\caption{a diagram which does not have an admissible quasi-Cartan companion} 

\label{fig:nonadm}

\end{figure}
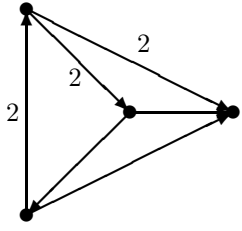

Let us also note that Theorem \ref{th:companion0} could be useful for recognizing mutation classes of acyclic diagrams: a diagram that does not have an admissible quasi-Cartan companion can not be obtained from any acyclic diagram by a sequence of mutations. An example of such a diagram is given in Figure~\ref{fig:nonadm}. (We refer to \cite[Section 2]{S3} for properties of diagrams of skew-symmetrizable matrices). Another application of the admissibility property to the corresponding Weyl groups can be found in \cite{S7}, where a fundamental class of relations have been shown to be satisfied by the reflections of the $\cc$-vectors.



\section{Proofs of main results}
\label{sec:proof}

Let us first recall the following well-known property of root systems: 
   For a generalized Cartan matrix $ A $ with symmetrizing matrix $ D={diag} (d_1,...,d_n) $,  there is an invariant symmetric bilinear form $ ( , ) $ defined on the simple roots as $ (\alpha_i,\alpha_j) =d_iA_{i,j}=d_jA_{j,i}=(\alpha_j,\alpha_i)$. 
  Let us note that, for any real root $ \alpha$, the corresponding reflection $ s_{\alpha} $ is defined on the real roots as
  $s_{\alpha}(\beta)=\beta-\langle\beta,{\alpha}^{\vee}\rangle\alpha$, with $\langle\beta,{\alpha}^{\vee}\rangle= \dfrac{2(\alpha,\beta)}{(\alpha,\alpha)} $. 
   In particular,  $s_{\alpha_i}(\alpha_j)=\alpha_j-\langle\alpha_j,{\alpha_i}^{\vee}\rangle\alpha_i=\alpha_j-A_{i,j}\alpha_i$. 
 
Let us also recall the mutation of quasi Cartan companions \cite[Definition 1.6]{S6}. Suppose that $B$ is a skew-symmetrizable matrix and let $A$ be a quasi-Cartan companion of $B$. 
Let $k$ be an index. For each sign $\epsilon=\pm1$, "the $\epsilon$-mutation of $A$ at $k$" is the quasi-Cartan matrix 
$\mu_{k}^\epsilon(A)=A'$ such that for any $i,j \ne k$: $A'_{i,k}=\epsilon sgn(B_{k,i})A_{i,k}$, $A'_{k,j}=\epsilon sgn(B_{k,j})A_{k,j}$, $A'_{i,j}=A_{i,j}-sgn(A_{i,k}A_{k,j})[B_{i,k}B_{k,j}]_+$. 
In the setup of Theorem \ref{th:DWZ2-RS}, suppose that $t \overunder{k}{} t'$  in $\TT_n$ and let $ A $ and $ A' $ be the associated quasi-Cartan companions. Then $A'=\mu^{\epsilon}_k(A)$ for $\epsilon=sgn(\cc_k)$.

We first prove Theorem~\ref{th:companion0} for convenience:

\noindent 
{Proof of Theorem~\ref{th:companion0}.}
To prove the first part, let us suppose that $sgn(B_{j,i})=sgn(\cc_j)$. 
Let $\mu_j(\cc,B)=(\cc',B')$ with $B'=\mu_j(B)$. Then $\cc_i'=\cc_i+[sgn(\cc_j)B_{j,i}]_+\cc_j=\cc_i+sgn(B_{j,i})B_{j,i}\cc_j=\cc_i+|B_{j,i}|\cc_j$. 
We denote by  $ ( , ) $ the invariant symmetric bilinear form defined by $ A_{0} $ on the root lattice and let $ D=diag (d_1,...,d_n) $ be the symmetrizing matrix for $ A_{0} $. 
Note that, by Theorem \ref{th:DWZ2-RS}, we have the following:
$2d_i=({\cc'_i},\cc'_i)=(\cc_{i},\cc_i)$, $2d_j=(\cc_{j},\cc_j)$, $(\cc_j,\cc_i)=(\cc_i,\cc_j)=d_iA_{{i,j}}=d_jA_{j,i}$.
Then
$2d_i=({\cc'_i},\cc'_i)=(\cc_i+|B_{j,i}|\cc_j,\cc_i+|B_{j,i}|\cc_j)=
(\cc_i,\cc_i)+(\cc_i,|B_{j,i}|\cc_j)+(|B_{j,i}|\cc_j,\cc_i)+|B_{j,i}|(\cc_j,|B_{j,i}|\cc_j)
=(\cc_i,\cc_i)+2|B_{j,i}|(\cc_j,\cc_i)+|{B}_{j,i}|^2(\cc_{j},\cc_j)
=2d_i+2|B_{j,i}|(\cc_{j},\cc_i)+|B_{j,i}|^22d_j=2d_i+2|B_{j,i}|d_jA_{j,i}+|B_{j,i}|^22d_j=2d_i+2|B_{j,i}|d_j(A_{j,i}+|B_{j,i}|)$, 
implying that $A_{j,i}+|B_{j,i}|=0$, thus $A_{j,i}=-|B_{j,i}|=-sgn (B_{j,i})B_{j,i}=-sgn (\cc_j)B_{j,i}$. 


To prove the second part of the theorem, let us suppose that $sgn(B_{j,i})=-sgn(\cc_j)$. 
Let $\mu_i(\cc,B)=(\cc',B')$ with $B'=\mu_i(B)$. Note that $sgn(B'_{j,i})=-sgn(B_{j,i})$ and $ |B'_{j,i}|=|B_{j,i}| $ (by the definition of mutation). Let $ A' $ be  the quasi-Cartan companion associated to the $ Y $-seed $(\cc',B')$ (Theorem \ref{th:DWZ2-RS} ), (Note then that $A'=\mu^{\epsilon}_i(A)$ where $\epsilon=sgn(\cc_i)$).

For the proof, we first assume that $sgn(\cc_j)=-sgn(\cc_i)$. Then we have $ sgn({\cc}_i)=sgn(B_{j,i}) $, so  ${\cc'}_j=\cc_j$ and ${\cc'}_i=-\cc_i$, implying $sgn({\cc'}_j)=sgn({\cc}_j)=-sgn(B_{j,i})=sgn(B'_{j,i})$, i.e. for the $Y$-seed $(\cc',B')$, we have $sgn(B'_{j,i})=sgn({\cc'}_j)$.
Thus, by the first part of the theorem, we have $-|B'_{j,i}|=A'_{j,i}=-A_{j,i}$. Thus $A_{j,i}=|B'_{j,i}|=|B_{j,i}|=-sgn(\cc_i)sgn(\cc_j)|B_{j,i}|$. 


Let us now assume that $sgn(\cc_j)=sgn(\cc_i)$. Then, since we have assumed $sgn(B_{j,i})=-sgn(\cc_j)$, we have $ sgn({\cc}_i)=-sgn(B_{j,i})= sgn(B_{i,j})$. Then, by the first part of the theorem, we have $A_{i,j}=-|B_{i,j}|$.  
Thus, since $ A $ is symmetrizable and a quasi-Cartan companion, we also have $A_{j,i}=-|B_{j,i}|$, which is equal to $-sgn(\cc_i)sgn(\cc_j)|B_{j,i}|$.


On the other hand, our assumption $sgn(B_{j,i})=-sgn(\cc_j)$ implies the following: $-sgn(\cc_i)sgn(\cc_j)|B_{j,i}|=-sgn(\cc_i)sgn(\cc_j)sgn(B_{j,i})B_{j,i}= \\
-sgn(\cc_i)sgn(\cc_j)(-sgn(\cc_j))B_{j,i}= sgn(\cc_i)B_{j,i}$. This completes the proof.

\vspace{0.1in}

\noindent 
{Proof of Corollary \ref{cor:companion}.}
Let us note that for $\mu_k(\cc,B)=(\cc',B')$ we have the following: $\cc'_k=-\cc_k$; $\cc'_i =\cc_i+[sgn(\cc_k)B_{k,i}]_+\cc_k$ if $i \neq k$ by \eqref{eq:y-mutation}. 
On the other hand, $[sgn(\cc_k)B_{k,i}]_+\ne 0$ if and only if $sgn(\cc_k)B_{k,i}>0$ if and only if $sgn(\cc_k)=sgn(B_{k,i})$. Then, by  Theorem~\ref{th:companion0}, we have $[sgn(\cc_k)B_{k,i}]_+=-A_{k,i}$. 
Thus $\cc'_i=\cc_i-A_{k,i}\cc_k=s_{\cc_k}(\cc_i)$ by the definition of a reflection. Also $\cc'_k=-\cc_k=s_{\cc_k}(\cc_k)$ This completes the proof of the statement.

\vspace{0.1in}

\noindent 
{Proof of Theorem~\ref{th:companion2}.} As we discussed in Section \ref{sec:intro}, the special case of this theorem when $ B $ is skew-symmetric was obtained in \cite[Theorem 1.4]{S6} by the author. The proof in \cite{S6} uses only the general properties of the mutations of skew-symmetrizable matrices with quasi-Cartan companions and the properties given in Theorem~\ref{th:companion0} (which was obtained for skew-symmetric matrices in \cite[Theorem 1.3]{S6}; note that in this case the companion $ A $ is symmetric and $A_{i,j}={\cc_{i}}^TA_0\cc_j$). 
Since we have proved Theorem~\ref{th:companion0} above for skew-symmetrizable matrices, the proof of \cite[Theorem 1.4]{S6} also holds for the skew-symmetrizable matrices. 
Thus, for the proof of Theorem~\ref{th:companion2}, we refer the reader to the proof of  \cite[Theorem 1.4]{S6}. 


\end{document}